\documentclass[preprint,12pt,authoryear]{elsarticle}

\usepackage{amssymb, amsmath}
\usepackage{amsthm}
\usepackage{hyperref}
\hypersetup{
    colorlinks=true,
    linkcolor=blue,
    }
\usepackage{stmaryrd}
\usepackage{subfig}
\graphicspath{ {./images/} }
\newcounter{example}[section]

\journal{arXiv}

\newcommand{\mV}{\ensuremath{{\mathbb V}}}
\newcommand{\mR}{\ensuremath{{\mathbb R}}}
\newcommand{\mZ}{\ensuremath{{\mathbb Z}}}

\newcommand{\cR}{\ensuremath{\mathcal R}}

\newcommand{\cd}{ \varobslash }

\newcommand{\B}{\ensuremath{\mathcal B}}

\newcommand{\C}{\ensuremath{\mathcal C}}
\newcommand{\A}{\ensuremath{\mathcal A}}

\newcommand{\G}{\ensuremath{\mathcal G}}

\renewcommand{\vec}[1]{\mbox{vec}(#1)}

\newdefinition{Ex}{Example}

\newcommand{\MARjP}{$MAR_j(p)$~}

\newcommand{\blue}[1]{\textcolor{blue}{#1}}
\newcommand{\red}[1]{\textcolor{red}{#1}}
\newcommand{\green}[1]{\textcolor{green}{#1}}

\begin{document}

\begin{frontmatter}


\title{On the relation between Global VAR Models and Matrix Time Series Models with Multiple Terms}


\author[BI]{Kurtulus Kidik}

\author[BI]{Dietmar Bauer}

\affiliation[BI]{organization={Econometrics, Bielefeld University},
            addressline={Universitätsstrasse 25}, 
            city={Bielefeld},
            postcode={33615}, 
            state={NRW},
            country={Germany}}

\begin{abstract}
Matrix valued time series (MaTS) and global vector autoregressive (GVAR) models both impose restrictions on the general VAR for multidimensional data sets, in order to bring down the number of parameters.

Both models are motivated from a different viewpoint such that on first sight they do not have much in common. When investigating the models more closely, however, one notices many connections between the two model sets. 

This paper investigates the relations between the restrictions imposed by the two models. 
We show that under appropriate restrictions in both models we obtain a joint framework allowing to gain insight into the nature of GVARs from the viewpoint of MaTS. 

\end{abstract}


\end{frontmatter}


\section{Introduction}
\label{sec:Intro}

In this paper, we consider a matrix times series (MaTS) $(Y_t)_{t \in \mZ}, Y_t \in \mR^{M \times N}$ where $M$ denotes the number of variables and $N$ the number of regions. Here, it is assumed that in all regions the same variables have been measured as would be typical for macro-economic data sets provided, for example, by the OECD \citep{OECD} or the FRED data base  of the Federal Reserve Bank of St. Louis \citep{FRED}. 
The $n$-th column $Y_{t,:,n} \in \mR^M$ then denotes the vector of observations for the region $n$ at the time point $t$. Throughout this paper, uppercase letters denote matrices (for example $Y_t$) and lowercase letters refer to the corresponding vectorized processes such that $(y_t)_{t \in \mZ}, y_t \in \mR^{MN}$, is the vector process obtained from columnwise vectorization of $Y_t$. 
We will also use the notation ':' to indicate rows or columns of a matrix such that, for example, $Y_{t,m,:}$ denotes the $m$-th row and $Y_{t,:,n}$ the $n$-th column of the matrix $Y_t$. 

In this notation, the matrix time series model with multiple terms (\MARjP) is formulated as: 

\begin{equation} \label{eq:MARjP}
Y_t = \sum_{j_1=1}^{J_1} A_{1,j_1}Y_{t-1}B_{1,j_1}' + \dots + \sum_{j_p=1}^{J_p} A_{1,j_p}Y_{t-p}B_{p,j_p}' + U_t, \quad t \in \mZ,
\end{equation} 

where $A_{i,j_i}\in \mR^{M \times M}, B_{i,j_i}\in \mR^{N \times N}$ and $(u_t)_{t \in \mZ}, u_t \in \mR^{MN}$ is a white noise process. 

Oftentimes the noise $U_t$ is assumed to follow a {\em separable} multivariate normal distribution in the sense that 
$\mV(u_t) = \Sigma_c \otimes \Sigma_r$ (see, for instance, \cite{Gupta2000}, \cite{CHEN2021539}). A leading case $\Sigma_c = I_N$ implies that the noise for 
each region is uncorrelated with noises from the other regions. The variance 
$\mV(u_t) = (I_N \otimes \Sigma_r)$ then is a block diagonal with equal diagonal 
blocks. 

Vectorizing this equation, we obtain a vector autoregressive (VAR) model: 

\begin{equation} \label{eq:vecMAR} 
y_t = \underbrace{\left( \sum_{j_1=1}^{J_1} (B_{1,j_1} \otimes A_{1,j_1}) \right)}_{\A_1} y_{t-1}+ \dots + \underbrace{\left( \sum_{j_p=1}^{J_p} (B_{p,j_p} \otimes A_{p,j_p}) \right)}_{\A_p}  y_{t-p}+ u_t, \quad t \in \mZ.
\end{equation} 

Noting that the number of elements in $\A_1 := \left( \sum_{j_1=1}^{J_1} (B_{1,j_1} \otimes A_{1,j_1}) \right)$ equals $(MN)^2$ while the number of elements in $A_{1,j},B_{1,j}$ equals $(M^2+N^2)J_1-J_1^2$ (taking identifiability restrictions into account, see below for details), $MAR_j(p)$ is a restriction of the VAR model, if $0 \le J_1 \le (\min(M,N))^2$. 

Thus, the $MAR_j(p)$ model is essentially a structured vector autoregression (VAR) specification. This note discusses the relations to other structured VAR specifications, in particular the GVAR model.

\section{The Global VAR model} \label{sec:GVAR} 
The global vector autoregressive model \citep{ChudikPesaranGVAR} introduces restrictions to the VAR in a different fashion: The main idea here is to decompose the model building process into separately specified models for each region:

\begin{equation} \label{eq:regional}
\begin{array}{lll}
Y_{t,:,i} & = & A_1^{(i)} Y_{t-1,:,i} + \dots + A_p^{(i)} Y_{t-p,:,i} \\
& & + 
 B_0^{(i)} Y_{t,:,i}^* + \dots + B_q^{(i)} Y_{t-q,:,i}^* 
 + U_{t,:,i} \quad i = 1, \dots, N
\end{array} 
\end{equation} 

where $A_j^{(i)} \in \mR^{M \times M}, B_j^{(i)} \in \mR^{M \times S}$ and 
$Y_{t}^* \in \mR^{S \times N}$ denotes exogenous variables (in an appropriate 
sense) representing spill-over effects.\footnote{The number of star variables may vary across regions, but for presentation reasons, we here use the same number, which may be achieved by padding the matrices with zeros to achieve equal numbers. Similarly, the number of variables within a region may also vary in general; however, in the MaTS environment, the same set of variables is modeled in each region.} 

Combining all regional models into a large vector, we obtain:
\begin{equation} \label{eq:GVAR} 
y_t = \A_{1,\cd} y_{t-1} + \cdots + \A_{p,\cd} y_{t-p} + 
\sum_{j=0}^q \B_{j,\cd} y_{t-j}^*  
+ u_t. 
\end{equation} 

Here,$\A_{j,\cd}= \mbox{diag}(A_j^{(1)},...,A_j^{(N)}) \in \mR^{NM \times NM}$ denotes a block diagonal matrix. 

The main idea behind GVARs is the definition of star variables $Y_{t,:,i}^*$ summarizing the influence from the remaining regions. These variables are typically defined as the weighted average of variables for neighboring regions:

\begin{equation} \label{eq:starvar}
Y_{t,m,n}^* = \sum_{j=1, j \ne n}^N w_{m,n,j} Y_{t,m,j}
= Y_{t,m,:} W_{m,n}, \quad w_{m,n,n}=0, W_{m,n} = \begin{pmatrix}
w_{m,n,1} \\ w_{m,n,2} \\ \vdots \\ w_{m,n,N}     
\end{pmatrix}\in \mR^N. 
\end{equation} 

For example, choosing $w_{m,n,j} = 1/(N-1), j \ne n,$ leads to averages over the variables of all other regions. If this is used for all variables ($w_{m,n,j} = \omega_{n,j}, W_{m,n}= W_n, \forall m$) in the system, we obtain: 

$$
Y_{t,:,n}^* = Y_{t} W_n.
$$

Again, stacking these definitions, we obtain: 

$$
y_t^* = W y_t = \begin{pmatrix}
0 & \omega_{1,2} I_M & \cdots & \omega_{1,N} I_M \\
\omega_{2,1} I_M & 0 & \ddots & \vdots \\
\vdots &  \ddots & \ddots & \vdots \\
\omega_{N,1} I_M & \cdots & \omega_{N,N-1}I_M & 0 
\end{pmatrix} y_t. 
$$

Note that $W = \tilde W \otimes I_M, \tilde W = [W_1,...,W_N]' \in \mR^{N \times N}$. If not all variables are represented in the star variables, $I_M$ can be replaced by $\tilde I_M$ of smaller dimension. 

A different specification uses a triangular dependence such that $\omega_{n,j}=0, j\ge n$, in which case $W$ is lower triangular. 
This specification introduces a special situation of influence, 
wherein regions are ordered such that regions earlier in the ordering influence the later regions, but not vice versa. 

Using the relation $y_t^* = Wy_t$ the model becomes (for simplicity of presentation assuming $p=q$ and under the assumption that $\G_0$ defined below is invertible)

\begin{align*}
y_t &= \sum_{j=1}^p \A_{j,\cd} y_{t-j}  + 
\sum_{j=0}^p \B_{j,\cd} W y_{t-j}  
+ u_t \\
&= \B_{0,\cd}W y_{t} + \sum_{j=1}^p (\A_{j,\cd} + \B_{j,\cd}W)y_{t-j}  + 
u_t \\
\Rightarrow \underbrace{(I_{MN} - \B_{0,\cd}W)}_{\G_0} y_t &= \sum_{j=1}^p (\A_{j,\cd} + \B_{j,\cd}W)y_{t-j} 
+ u_t, \\
y_t &=  \sum_{j=1}^p \G_0^{-1}(\A_{j,\cd} + \B_{j,\cd}W)y_{t-j}  
+ \G_0^{-1}u_t.
\end{align*}

It follows that we obtain a VAR. If $\B_{0,\cd} \ne 0$ it is given in what we could call 'structural form'. For a triangular system, $\G_0$ is guaranteed to be invertible, while otherwise invertibility is not trivial and must be assumed (see \cite{ChudikPesaranGVAR}). 

For a triangular system, where $W$ is block lower triangular, all matrices in the vectorized equation are block lower-triangular. 

\section{Properties of the \MARjP} \label{sec:MArJP} 

The \MARjP requires additional restrictions in order to be identified. 
Clearly, the Kronecker product $B \otimes A$ identifies the two matrices $A$ and $B$ 
only up to a scalar multiple. When more than one term is involved (i.e. $J_r >1$), additional issues arise. In order to understand the structure of \MARjP models, one can use two main tools, that are described below. 

{\bf Tool 1: the magical mapping:} A key tool in ascertaining an identifiable parameterization is given by the mapping $\cR: \mR^{MN \times MN} \to \mR^{M^2 \times N^2}$ (sometimes called 'rearrangment operator' \cite{VanLoan1993}) defined via:

$$
\cR ( B \otimes A ) = \vec{A} \vec{B}'. 
$$

Vectorising the matrices, one can verify that the mapping rearranges the entries in the 
matrices. This can be demonstrated for $M=N=2:$

\begin{align*}
A = \begin{pmatrix}
    a_{11} & a_{12} \\ a_{21} & a_{22} 
\end{pmatrix}, &B =     \begin{pmatrix}
    b_{11} & b_{12} \\ b_{21} & b_{22} 
\end{pmatrix}, \\
B \otimes A &= \left( \begin{array}{cc|cc}
\blue{a_{11} b_{11}} & \blue{a_{12} b_{11}} & \red{a_{11} b_{12}} & \red{a_{12} b_{12}} \\
\blue{a_{21} b_{11}} & \blue{a_{12} b_{11}} & \red{a_{11} b_{12}} & \red{a_{12} b_{12}} \\ \hline 
\green{a_{11} b_{21}} & \green{a_{12} b_{21}} & a_{11} b_{22} & a_{12} b_{22} \\
\green{a_{21} b_{21}} & \green{a_{22} b_{21}} & a_{21} b_{22} & a_{22} b_{22}
\end{array} \right), \\
\cR ( B \otimes A) &= \left(
\begin{array}{c|c|c|c}
\blue{a_{11} b_{11}} & \green{a_{11} b_{21}} & \red{a_{11} b_{12}} & a_{11} b_{22} \\
\blue{a_{21} b_{11}} & \green{a_{21} b_{21}} & \red{a_{21} b_{12}} & a_{21} b_{22} \\
\blue{a_{12} b_{11}} & \green{a_{12} b_{21}} & \red{a_{12} b_{12}} & a_{12} b_{22} \\
\blue{a_{22} b_{11}} & \green{a_{22} b_{21}} & \red{a_{22} b_{12}} & a_{22} b_{22}
\end{array} \right).
\end{align*}

The mapping, therefore, is bijective and linear (it can be described by a $(MN^2) \times (MN)^2$ selector matrix of zeros and ones linking the vectorizations) such that 

$$
\cR \left(\sum_{j=1}^J  (B_j \otimes A_j ) \right) = \sum_{j=1}^J \vec{A_j}\vec{B_j}'. 
$$

Note that each $\vec{A_j}\vec{B_j}'$ is a matrix of rank-1 such that in general the sum of $J$ matrices leads to a matrix of rank-$J$. 

This representation can be used in order to achieve identifiability: 
\cite{li2021multilinearTenAr} and \cite{hsu2024rank} provide identifiability conditions based on the singular value decomposition (SVD) of the rearranged matrix 
$\cR(\A_1)$. Alternatively, one may use the QR decomposition to identify the factors $\A$ and $\B$ from the product $\cR(\A_1) = \A \B'$. Here $\A = [ \vec{A_1},\vec{A_2},...,\vec{A_J}] \in \mR^{M^2 \times J}, 
\B = [ \vec{B_1},\vec{B_2},...,\vec{B_J}] \in \mR^{N^2 \times J}$ such that 
$\A'\A = I_{J}$ and $\B'$ is positive upper triangular (see \cite{bauer2020parameterization}, Lemma 1). This approach is followed in \cite{Kidik}.

{\bf Tool 2: blockwise decomposition:} An alternative insight can be obtained from the blockwise representation:

\begin{align*}
B \otimes A &= \begin{pmatrix}
b_{11} A & b_{12} A & \cdots & b_{1N} A \\ 
b_{21} A & \ddots & \ddots & \vdots \\
\vdots & \ddots & \ddots & \vdots \\
b_{N1} A & \cdots & \cdots & b_{NN} A 
\end{pmatrix} \Rightarrow \\~\\
\cR ( B \otimes A ) &= \begin{pmatrix}
  b_{11} \vec{A} & b_{21} \vec{A} & \cdots & b_{NN} \vec{A}  
\end{pmatrix}, \\~\\
\C &= \begin{pmatrix}
C_{11} & C_{12}  & \cdots & C_{1N}  \\ 
C_{21} & \ddots & \ddots & \vdots \\
\vdots & \ddots & \ddots & \vdots \\
C_{N1} & \cdots & \cdots & C_{NN}      
\end{pmatrix}
\Rightarrow \\~\\
\cR ( \C ) &= \begin{pmatrix}
  \vec{C_{11}} &  \vec{C_{21}} & \cdots & \vec{C_{NN}}  
\end{pmatrix}
\end{align*}

This implies, for example, that every lower triangular matrix can be written as the sum of at most $N(N+1)/2$ Kronecker product terms. $\G_0$ in the lower triangular case consists of a maximum of $N(N-1)/2+1$ 
terms, since the diagonal blocks are equal to unity. 

Moreover, in the case $W = \tilde W \otimes \tilde I_M$ we obtain the following. 

$$
\B_{0,\cd} W = \begin{pmatrix} \omega_{11} B_{0}^{(1)} \tilde I_M & \cdots & \omega_{1N} B_{0}^{(1)} \tilde I_M \\ 
\vdots & & \vdots \\ 
\omega_{N1} B_{0}^{(N)} \tilde I_M & \cdots & \omega_{NN} B_{0}^{(N)} \tilde I_M
\end{pmatrix} 
$$

This matrix can be written as the sum of at most $N$ Kronecker product terms so that $\G_0$ can be written as the sum of at most $N+1$ terms. 

\section{Links between GVAR and \MARjP}
As both GVAR and \MARjP are restrictions of the general VAR system for the vectorized 
matrix time series and the \MARjP for the maximal number of terms is equivalent to 
the VAR, clearly the GVAR can be obtained as a special case of a \MARjP. In particular, when the star variables are defined using $W = \tilde W \otimes \tilde I_M$ we get from above

$$
\G_0 y_t = \sum_{j=1}^p (\A_{j,\cd} + \B_{j,\cd}W)y_{t-j} 
+u_t
$$

where 

$$
\A_{j,\cd} + \B_{j,\cd}W = \begin{pmatrix}
A_j^{(1)} & \omega_{1,2} B_j^{(1)}  &   \cdots & \omega_{1,N} B_j^{(1)}  \\  
\omega_{2,1} B_j^{(2)} & A_j^{(2)} & \cdots & \vdots \\
\vdots & \cdots & \ddots & \vdots \\
\omega_{N,1} B_j^{(N)} & \cdots & \omega_{N,N-1} B_j^{(N)} & A_j^{(N)} 
\end{pmatrix}. 
$$

is a sum of Kronecker products with a maximum of $2N$ terms: a maximum of $N$ for the block diagonal matrix $\A_{j,\cd}$ and an additional maximal $N$ for $\B_{j,\cd}W$. This is typically much smaller than $\min(M^2,N^2)$. 

The inclusion of $\G_0$ on the left hand side is non-standard for the \MARjP 
models. Multiplying with $\G_0^{-1}$, one obtains a reduced form equation. This can be estimated using the \MARjP framework. One obtains the following equation: 

$$
y_t = \sum_{j=1}^p \hat \A_j y_{t-j} + \hat u_t
$$

where $\hat \A_j = \sum_{r=1}^{J_j} (\hat B_{j,r} \otimes \hat A_{j,r})$. From these reduced form estimates one obtains equations to identify the structural form parameters:

$$
\hat \G_0 \hat \A_j = \hat\A_{j,\cd} + \hat\B_{j,\cd}W + \hat V_j, \quad j=1,...,p. 
$$

The uncorrelatedness of the star variables with the residuals of the $i$-th region adds one block of equations: 

$$
\sum_{j=1, j \ne i}^N w_{i,j} B_0^{(i)} \sum_{t=p+1}^T Y_{t,:,j} \hat U_{t,:,i}' = 0.
$$

Counting the number of restrictions, we obtain $p(MN)^2 + NM^2$ equations for 
the $2pNM^2$ (for $A_j^{(i)}, B_j^{(i)}$) plus $NM^2$ (for $B_0^{(i)}$) plus $N(N-1)$ (for $w_{i,j}$) parameters. The difference therefore equals

$$
(N-2)pNM^2 - N(N-1) = N [(N-2)pM^2 - N+1 ]
$$

which is in general positive for $N>2$. Therefore, from these equations one can obtain
a least absolute distance or a least squares estimator for all parameters. Restrictions 
$\sum_{j=1, j\ne i}^N w_{i,j} = 1$, $w_{i,i} = 0 $ and $w_{i,j} \ge 0 $ can be imposed. 

Starting from this initial estimator, one may use an alternating optimization routine to estimate the GVAR parameters including the weights sequence:

\begin{itemize}
    \item Estimate $A_j^{(i)}, B_j^{(i)}$ given $W$: This is a standard regional GVAR estimate.
    \item Estimate $W$ given $A_j^{(i)}, B_j^{(i)}$: 
\begin{align*}
Y_{t,:,i} - \sum_{j=1}^p A_j^{(i)} Y_{t-j,:,i} &= 
 B_0^{(i)} Y_{t,:,i}^* + \dots + B_q^{(i)} Y_{t-q,:,i}^* 
 + U_{t,:,i} \\
 &= \sum_{j=1, j \ne i}^N w_{i,j} ( B_0^{(i)} Y_{t,:,j} + \dots + B_q^{(i)} Y_{t-q,:,j} ) + U_{t,:,i} \\
 &= \sum_{j=1, j \ne i}^N w_{i,j} \tilde Y_{t,:,j}  + U_{t,:,i} 
\end{align*} 
\end{itemize}

In the second equation, the error terms have variance $\Omega_i$ which can be taken into 
account. Also, the regressors $\tilde Y_{t,:,j}$ can be endogenous, hence IV estimation 
should be used. 

Moreover, for the error term, we obtain for the reduced form errors $\G_0^{-1} u_t$ such that the error variance equals (if the matrix $\mV(u_t)=\Sigma_{\cd} $ is block diagonal)

$$
\G_0^{-1} \mV(u_t) \G_0^{-T} = \G_0^{-1} \Sigma_{\cd} \G_0^{-T}
$$

In the triangular case, we hence get that the Cholesky factor of the innovation 
variance equals $\G_0^{-1} \Sigma_{\cd}^{1/2}$. Since in $\G_0$ the block diagonal 
equals the identity matrix, it is simple to identify $\G_0^{-1}$ and 
$\Sigma_{\cd}^{1/2}$ from the Cholesky factor in that case. 


While the same procedure could be followed starting from VAR estimates $\hat \A_j$ the number of free entries in the VAR is much larger than in the approximating \MARjP. As such that \MARjP can be seen as providing an intermediate model between the highly restricted GVAR and the unrestricted VAR that allows to obtain more precise estimates to start the GVAR estimation (including estimation of the weights). Additionally, one can  
test the appropriateness of the GVAR within the smaller model and hence with larger power compared to tests using the unrestricted VAR. 
%
%
%
%

\section{Conclusions} \label{sec:concl}

This note shows that the \MARjP with a limited number of terms can be used as an intermediate model between the heavily restricted GVAR and the unrestricted VAR model. The latter is hard to estimate for typical sample sizes due to the large number of parameters involved, while the former may be too strongly restricted. In addition, the weightings $\omega_{i,j}$ involved in the GVAR model are often not known and the definition involves heuristic arguments. The \MARjP provides some middle ground between these two extremes, providing help with the specification of the GVAR and in particular ways to estimate the weights. 



\bibliographystyle{plainnat}
\bibliography{lit_econometrics}

\end{document}